\documentclass[reqno]{amsart}
\usepackage{hyperref}

\begin{document}
\title[
\hfilneg %EJDE-2005/XX
\hfil Bounded solutions with negative stiffness]
{On bounded solutions for second order linear differential equations with negative stiffness}
\author[C\'esar A. Terrero-Escalante
\hfil %EJDE-2005/XX
\hfilneg] 
{C\'esar A. Terrero-Escalante}  % in alphabetical order

\address{Departmento de F\'{\i}sica Te\'orica\\ 
Instituto de F\'{\i}sica \\ 
Universidade do Estado do Rio de Janeiro\\ 
Rua S\~ao Francisco Xavier 524,
Maracan\~a, 20559-900 Rio de Janeiro, RJ, Brazil.}
\email[C\'esar A. Terrero-Escalante]{cterrero@dft.if.uerj.br}

\date{}
%\thanks{Submitted November XX, 2005. Published XXXX X, 200X.}
\begin{abstract}
 Second order linear non-autonomous differential equations with negative stiffness are considered. Using
 Chetaev-like (Lyapunov-like) functions, necessary (sufficient) conditions are found for the
 solutions to be bounded for all initial conditions if any one of the coefficients is constant. 
 Conclusions are then extended to include systems where both coefficients are time-varying. 
\end{abstract}

\maketitle
\numberwithin{equation}{section}
\newtheorem{theorem}{Theorem}[section]
\newtheorem{lemma}[theorem]{Lemma}
\newtheorem{corollary}[theorem]{Corollary}
\newtheorem{remark}[theorem]{Remark}
\newtheorem{example}[theorem]{Example}

\section{Introduction}

The second order linear differential non-autonomous homogeneous equation  
\begin{equation}
\ddot{x}+\alpha(t)\dot{x}-k(t)x=0 \,, 
\label{eq:Main}  
\end{equation} 
and its non-homogeneous version are usually found in mechanical problems like the study of vibrations,  
and in perturbative analysis of solutions to differential equations in Physics.  
Here, $\alpha,k:\textsl{I}\rightarrow\mathbb{R}$,  
dot stands for differentiation with respect to the independent variable  
$t \in \textsl{I}\equiv [t_0,\infty)$.  
Following Mechanics notation,  
$\alpha(t)$ is going to be called the {\it damping} of the system,  
while $-k(t)$ will be referred as the {\it stiffness}.  

Too often a $k(t)$ strictly positive for all $t\geq t_0$ is  
automatically and tacitly  
regarded as a signal of instability of the trivial solution of the above equation.  
The reason for that is just that the corresponding to Eq.(\ref{eq:Main}) planar system,  
\begin{equation} 
\dot{\bf x}={\bf x}^T 
\,\begin{array}{|l|l} 0 \quad \hspace{1.cm} 1 \\ \\ k(t) \quad -\alpha(t) \end{array}\,\, 
{\bf x}  
\label{eq:MainSystem}  
\end{equation} 
(${\bf x}^T =| x_1 x_2 |$, $x_1\equiv x$, $x_2 \equiv \dot{x}$), 
has eigenvalues $\lambda_\pm(t) \propto -\alpha(t) \pm \sqrt{\alpha(t)^2+4k(t)}$,  
pointing to a saddle-like instability for any given value of $t$.  
Thus, if these eigenvalues are slowly changing the system is assumed to have no bounded (finite) solutions.  

However, a rigorous foundation of this folk theorem is still missing  
as well as a precise definition of {\it slowly changing}.  
A number of authors has addressed this problem  
(see for instance \cite{Rosenbrock,Desoer,Solo,Ignatiev,Juan})  
but either they deal only with positive stiffness  
or their results are not expressed as direct conditions on the time-varying coefficients.  

The aim of this paper is  
to help to fill this gap  
in the analysis of the stability on non-autonomous linear systems.  
With that in mind,  
in the next section,  
the Lyapunov's direct method is used to state  
two theorems for the necessary and sufficient conditions  
on the time variation of $\alpha(t)$ and $k(t)$,  
for Eq.(\ref{eq:Main}) to have bounded solutions for any initial condition ${\bf x}_0\equiv {\bf x}(t_0)$.  
Each theorem corresponds to a case where one of the coefficients is time-varying while the other one remains constant. 
In the last section,  
corollaries and examples are presented  
allowing to extend the conclusions to cases where both coefficients are time-varying.  
Throughout this manuscript,  
it is assumed that conditions are imposed on (\ref{eq:Main})  
such that existence and uniqueness of its solutions are warranted for all initial conditions  
${\bf x_0}\in \mathbb{R}^2$ and for all initial times $t_0$.  
The Euclidean norm, $\|\cdot\|$, in $\mathbb{R}^2$ is used.  
The definitions and theorems for boundedness and stability in the sense of Lyapunov  
can be found in many textbooks  
(see for instances \cite{Arnold,Khalil}). 
 
\section{Main results\label{sec:Main}} 
\begin{theorem} 
Consider the linear non-autonomous equation 
\begin{equation} 
\ddot{x}+\beta\dot{x}-k(t)x=0 \,, 
\label{eq:beta}  
\end{equation} 
where $\beta$ is a real constant and $k(t)$ is a differentiable real-valued function strictly positive in $\textsl{I}$. Then,  
\begin{enumerate} 
	\item $\beta>0$, \label{cond:Main1} 
	\item $\dot{k}(t)<0$ \label{cond:Main2} 
\end{enumerate} 
are necessary conditions for the solution ${\bf x}=0$ of (\ref{eq:beta}) to be stable in the sense of Lyapunov. Moreover,  
\begin{enumerate} 
\addtocounter{enumi}{2} 
		\item $\dot{\left[\frac 1 {k(t)}\right]}>\frac 2 {\beta}$\,, \label{cond:Main3} 
\end{enumerate} 
along with the conditions above, is sufficient for all solutions to be bounded $\forall {\bf x}_0$. 
\label{theor:Main} 
\end{theorem} 
\begin{proof} 
Condition (\ref{cond:Main1}) follows straightforwardly from  
requiring no continuous increase of the volume of the phase flow due to the linear vector field ${\bf v}=A(t){\bf x}$. 
In other words \cite{Arnold},  
${\rm div} {\bf v}={\rm Tr} A(t) \leq 0$,  
where ${\rm div}$ stands for the vector field divergence and ${\rm Tr}$ for the matrix trace.  
If $\beta=0$ the solutions will be unbounded  
since even in the limit $k(t)\rightarrow +0$ the solution blows up linearly.  
 
To prove condition (\ref{cond:Main2})  
recall first that for $\dot{k}(t)= 0$  
the eigenvalues of the corresponding planar system  
are $\lambda_\pm \propto -\beta \pm \sqrt{\beta^2+4k}$,  
and ${\bf x}=0$ is indeed a saddle point.  
So, assume now $\dot{k}(t)> 0$.  
Then, consider the Chetaev's function \cite{Khalil}, 
\begin{equation} 
V_\chi({\bf x},t)=\frac 1 2 k(t)x_1^2 - \frac 1 2 x_2^2 \,. 
\label{eq:Vchi} 
\end{equation} 
Let $D_1$ be the circle defined by $\|{\bf x}\|=\epsilon^2$,  
let $D_2$ be the union of the two subspaces distributed along $x_2=0$  
with borders $\partial D_2$ given by $x_2=\pm\sqrt{k(t_0)}x_1$,  
and let $D\subseteq D1\cap D2/\partial D_2$.  
Then, for $k(t)>0$, and $\epsilon$ as small as desired,  
$\forall {\bf x}\in D$: $V_\chi({\bf x},t)\geq V_\chi({\bf x},t_0)> 0$. 
$V_\chi({\bf x},t)=0$ only for ${\bf x}\in \partial D2$,  
what includes the origin.  
The time derivative of (\ref{eq:Vchi}) along the direction of the vector field ${\bf v}$ is 
\begin{equation} 
\frac {dV_\chi}{dt}=\frac 1 2 \dot{k}(t)x_1^2 + \beta x_2^2 \,. 
\label{eq:VchiDot} 
\end{equation} 
For any ${\bf x}\in D$ and for all $t\in \textsl{I}$: $\dot{V}_{\chi}(0,t)=0$ and $\dot{V}_{\chi}({\bf x},t)>0$.  
Therefore, any solution of (\ref{eq:beta}) with $\beta>0$, $\dot{k}(t)> 0$ and initial condition ${\bf x}_0\in D$ 
will eventually escape to infinity, independently of how close ${\bf x}_0$ is chosen to be to the origin. 
With these conditions, the trivial solution of (\ref{eq:beta}) is unstable.  
 
The sufficient condition (\ref{cond:Main3}) can be proved using the Lyapunov's function \cite{Arnold,Khalil} 
\begin{equation} 
V_\Lambda({\bf x},t)=\frac 1 2 k(t)x_1^2 + \frac 1 2 x_2^2 \,. 
\label{eq:VL} 
\end{equation} 
For $k(t)>0$ this function is non-negative and is equal to zero only at the origin.  
Since $k(t)$ is differentiable in $[t_0,\infty)$ and decreases with time,  
$\forall {\bf x}\neq 0$ and $\forall t\in \textsl{I}$ : $V_\Lambda({\bf x},t_0)> V_\Lambda({\bf x},t)>0$. 
$V_\Lambda({\bf x},t_0)\rightarrow\infty$ only when $\|{\bf x}\|\rightarrow\infty$.  
The corresponding time derivative along the direction of the vector field ${\bf v}$ is, 
\begin{equation} 
\frac {dV_\Lambda}{dt}=\frac 1 2 \dot{k}(t)x_1^2 + 2k(t)x_1x_2 - \beta x_2^2 \,, 
\label{eq:VLDot} 
\end{equation} 
which can also be written as, 
\begin{equation} 
\frac {dV_\Lambda}{dt}={\bf x^T} 
\, \begin{array}{|l|l} \frac 1 2 \dot{k}(t) \quad k(t) \\ \\ k(t) \quad -\beta  
\end{array}\,\, 
{\bf x} \,. 
\label{eq:VLDotM} 
\end{equation} 
This function is negative definite if  
$\forall t\geq t_0$: $\dot{k}(t)<0$, $\beta>0$ and $-\dot{k}(t)/k^2(t)>2 /\beta$.  
With these conditions, solutions of system (\ref{eq:beta}) will be bounded $\forall {\bf x}_0$. 
\end{proof} 
 
\begin{theorem} 
Consider the linear non-autonomous equation  
\begin{equation} 
\ddot{x}+\alpha(t)\dot{x}-\omega^2x=0 \,, 
\label{eq:omega}  
\end{equation} 
where $\omega$ is a real constant and $\alpha(t)$ is a differentiable real-valued function in $\textsl{I}$.  
Then,  
\begin{enumerate} 
    \item $\alpha(t)>0$, \label{cond:Main21} 
		\item $\dot{\alpha}(t)>0$, \label{cond:Main22} 
		\item $\frac d{dt} \ln \alpha(t)> \frac 1 2 \left[\frac 1 {\alpha(t)} + \omega^2 \right]^2$ \label{cond:Main23} 
\end{enumerate} 
$\forall t\geq t_0$ are together sufficient conditions  
for solutions of (\ref{eq:omega}) to be bounded $\forall {\bf x}_0$. 
\label{theor:Main2} 
\end{theorem} 
\begin{proof} 
This theorem can be proved using the Lyapunov function, 
\begin{equation} 
V_\Lambda({\bf x},t)=\frac 1 2 \frac 1 {\alpha(t)}x_1^2 + \frac 1 2 x_2^2 \,. 
\label{eq:VL2} 
\end{equation} 
and the corresponding time derivative along the direction of the vector field ${\bf v}$,  
\begin{equation} 
\frac {dV_\Lambda}{dt}={\bf x^T}  
 \,\begin{array}{|l|l} \frac 1 2 \frac d {dt}\left[\frac 1 {\alpha(t)}\right] \quad \frac 1 2 \left[\frac 1 {\alpha(t)}+\omega^2\right] \\ \\ \frac 1 2 \left[\frac 1 {\alpha(t)}+\omega^2\right] \quad -\alpha(t)  
\end{array}\,\, 
{\bf x} \,. 
\label{eq:L2} 
\end{equation} 
\end{proof} 
 
\section{Corollaries and examples} 
 
\begin{corollary} 
There exist time-varying coefficients $\alpha(\tau)$ and $k(\tau)$  
such that solutions to Eq.(\ref{eq:Main}) will be bounded $\forall {\bf x}_0$. 
\label{cor:first} 
\end{corollary} 
\begin{proof} 
\begin{example} 
The non-autonomous equation  
\begin{equation} 
x^{\prime\prime}+\alpha(\tau)x^{\prime}-k(\tau)x=0 \,, 
\label{eq:Ex1}  
\end{equation} 
where prime stands for derivative with respect to $\tau$,  
and $\forall\tau\geq \tau_0$ 
\begin{eqnarray} 
k(\tau)&=&\omega^2 f(\tau)^{n+1}\,,\\ 
\alpha(\tau)&=&\beta f(\tau)^{\frac n 2}-\frac n 2 f(\tau)^{-1}f(\tau)^{\prime}\,, 
\end{eqnarray} 
for all real $\omega$, $n>0$, $\beta>0$, $f(\tau)>0$, $f(\tau)^{\prime}<0$, and  
$f(\tau)^{-\frac n 2}[1/f(\tau)]^{\prime}>2 \omega^2/\beta$.  
\end{example} 
To test this, change the independent variable in Eq.(\ref{eq:beta}) using $d\tau/dt=f^{-n/2}$.  
\end{proof} 
\begin{corollary} 
In Eq.(\ref{eq:Main}) the time-varying coefficient $\alpha$ does not need to be strictly positive  
for the solutions to be bounded $\forall {\bf x}_0$. 
\label{cor:scnd} 
\end{corollary} 
\begin{proof} 
\begin{example} 
The non-autonomous equation  
\begin{equation} 
x^{\prime\prime}+\alpha(\tau)x^{\prime}-k(\tau)x=0 \,, 
\label{eq:Ex2}  
\end{equation} 
where $\forall\tau\geq \tau_0$, and for all real $\omega$, $n>0$, $\beta>0$, 
\begin{eqnarray} 
k(\tau)&=&\omega^2 g(\tau)^{n}\,,\\ 
\alpha(\tau)&=&\beta g(\tau)^{\frac n 2 + 1}-\frac n 2 g(\tau)^{-1}g(\tau)^{\prime}\,, 
\end{eqnarray} 
if $g(\tau)>0$, $g(\tau)^{\prime}>0$, and  
$d\ln g(\tau)/d\tau> g(\tau)^{n/2} [\beta^{-1} g(\tau)^{-1} + \omega^2 ]^2/2$.  
\end{example} 
And this can be tested using the parametric change $d\tau/dt=g^{-n/2}$ in Eq.(\ref{eq:omega}). 
\end{proof} 
Intuitively, a positive damping can control the instability induced by a negative stiffness 
if the former is asymptotically larger than the absolute value of the last. 
Theorems \ref{theor:Main} and \ref{theor:Main2}, together with corollaries \ref{cor:first} and \ref{cor:scnd} 
give a notion of how large the ratio $|\alpha(t)|/k(t)$ must be at any given time 
for the corresponding solution to be bounded $\forall {\bf x}_0$.  
\begin{remark} 
The sufficient condition \ref{theor:Main}.\ref{cond:Main3}  
for the existence of bounded solutions of Eq.(\ref{eq:beta}) 
is rather weak. 
\end{remark} 
Indeed, if the stiffness decays with respect to some finite positive damping $\beta$  
the solution of Eq.(\ref{eq:beta}) will asymptotically look like $C_1 + C_2\exp(-\beta t/2)$.  
For instance, for equation, 
\begin{example} 
\begin{equation} 
\ddot{x}+\beta\dot{x}-\omega^2{\rm t}^{-m}x=0 \,  
\label{eq:Ex3}  
\end{equation} 
with $m\geq 1$, condition \ref{theor:Main}.\ref{cond:Main3} reads 
\[ 
t^{m-1}>\frac 2 m \frac {\omega^2}{\beta}\, . 
\] 
So, the solutions are bounded even for $m=1$ if $\beta>2\omega^2$. 
For $m\geq 2$ the condition is satisfied asymptotically for any $\beta$ and $\omega$. 
\end{example} 
For $m=2$ the corresponding solutions are given by 
\[ 
x(t)=C_1 \sqrt{t} {\rm e}^{-\frac{\beta}2 t}\mathcal{I}\left(\frac{\sqrt{1+4\omega^2}}2,\frac{\beta t}{2}\right)  
+ C_2 \sqrt{t} {\rm e}^{-\frac{\beta}2 t}\mathcal{K}\left(\frac{\sqrt{1+4\omega^2}}2,\frac{\beta t}{2}\right)\, , 
\] 
where $C_1$ and $C_2$ are the integration constants, and $\mathcal{I}(a,z)$ and $\mathcal{K}(a,z)$  
are the modified Bessel functions of the first and second kinds, respectively.  
In the asymptotic regime, this solution can be approximated by 
\begin{equation} 
C_1 \frac 1 {\sqrt{\pi\beta}} + C_2 \frac {\sqrt{\pi}}{\sqrt{\beta}} {\rm e}^{-\beta t}\,, 
\end{equation} 
converging uniformly with $t\rightarrow\infty$ to the constant given by the initial condition $x_0$. 
For completeness, note that  
in the limit of strong decay of the negative stiffness, 
\begin{example} 
\begin{equation} 
\ddot{x}+\beta\dot{x}-\omega^2{\rm e}^{-mt}x=0 \, . 
\label{eq:Ex4}  
\end{equation} 
\end{example} 
the solutions for $m>0$ are given by, 
\begin{equation} 
x(t)=C_1 {\rm e}^{-\frac{\beta}2 t}\mathcal{I}\left(-\frac{\beta}{m},\frac{2\omega{\rm e}^{-\frac{m}2 t}}{m}\right) 
+ C_2 {\rm e}^{-\frac{\beta}2 t}\mathcal{K}\left(\frac{\beta}{m},\frac{2\omega{\rm e}^{-\frac{m}2 t}}{m}\right) 
\nonumber \,, 
\label{eq:SolEx4}  
\end{equation} 
which for very large $t$ can be approximated by 
\begin{equation} 
\frac 1 2 C_1 \left(\frac \beta m - 1\right)! \left(\frac \omega m \right)^{\left(-\frac \beta m\right)} 
+ C_2 \left[\left(\frac \beta m \right)!  
\left(\frac \omega m \right)^{\left(-\frac \beta m\right)}\right]^{-1}{\rm e}^{-\beta t}\,, 
\nonumber 
\label{eq:SolAsympEx4}  
\end{equation} 
converging uniformly with $t\rightarrow\infty$ to the constant given by $x_0$.  
\begin{remark} 
Conversely, the sufficient condition \ref{theor:Main2}.\ref{cond:Main23}  
for the existence of bounded solutions of Eq.(\ref{eq:omega}) 
is very strong, it is close to necessary. 
\end{remark} 
For instance, for equation, 
\begin{example} 
\begin{equation} 
\ddot{x}+\beta{\rm t}^{m}\dot{x}-\omega^2x=0\, , 
\label{eq:Ex5}  
\end{equation} 
with $m> 0$, condition \ref{theor:Main2}.\ref{cond:Main23} reads 
\[ 
t^{-1}-\frac 2 m \left[\frac {{\rm t}^{-m}}{\beta} + \omega^2 \right]^2  >0\, , 
\] 
which is not satisfied for any real $\beta$ and $\omega$ if $t\rightarrow\infty$. 
\end{example} 
In fact, for $0<m\leq 3$, Eq.(\ref{eq:Ex5}) can be solved in term of Heun series,  
each one diverging for all positive $\beta$ and real $\omega$. 
Nevertheless, this condition is not actually necessary  
since in the limit of strong growth for the damping, 
equation 
\begin{example} 
\begin{equation} 
\ddot{x}+\beta{\rm e}^{mt}\dot{x}-\omega^2x=0 \,, 
\label{eq:Ex5x}  
\end{equation} 
\end{example} 
has solution, 
\begin{eqnarray} 
x(t)&=& C_1 {\rm e}^{\frac 1 2 \left(mt - \frac \beta m {\rm e}^{mt}\right)} 
\left[\mathcal{I}\left(-\frac{1}{2}+\frac{\omega}{m},\frac 1 2 \frac \beta m {\rm e}^{mt}\right)  
+ \mathcal{I}\left(\frac{1}{2}+\frac{\omega}{m},\frac 1 2 \frac \beta m {\rm e}^{mt}\right)\right] \nonumber\\  
&+& C_2 {\rm e}^{\frac 1 2 \left(mt - \frac \beta m {\rm e}^{mt}\right)} 
\left[\mathcal{K}\left(\frac{1}{2}-\frac{\omega}{m},\frac 1 2 \frac \beta m {\rm e}^{mt}\right)  
- \mathcal{K}\left(\frac{1}{2}+\frac{\omega}{m},\frac 1 2 \frac \beta m {\rm e}^{mt}\right)\right] 
\nonumber \,, 
\label{eq:SolEx5}  
\end{eqnarray} 
for any real $\omega$, $\beta>0$ and $m>0$.  
Asymtotically it behaves like 
\begin{equation} 
C_1 \frac 1 {\sqrt{\pi}} \sqrt{\frac m \beta}  
+ C_2 {\sqrt{\pi}} \sqrt{\frac m \beta} {\rm e}^{- \frac \beta m {\rm e}^{mt}} \,, 
\label{eq:SolAsympEx5}  
\end{equation} 
also converging uniformly with $t\rightarrow\infty$ to the constant related to $x_0$. 
 
Finally, note that substitution $d\tau/dt=\tau$ (with $\tau_0>0$) will change Eq.(\ref{eq:Ex4}) into 
\begin{example} 
\begin{equation} 
x^{\prime\prime}+(\beta+1)\frac 1 \tau x^{\prime}-\omega^2\frac 1 {\tau^{m+2}}x=0 \,, 
\label{eq:Ex6}  
\end{equation} 
\end{example} 
and Eq.(\ref{eq:Ex5}) respectively into 
\begin{example} 
\begin{equation} 
x^{\prime\prime}+(\beta\tau^{m}+1)\frac 1 \tau x^{\prime}-\omega^2\frac 1 {\tau^{2}}x=0 \,. 
\label{eq:Ex7}  
\end{equation} 
\end{example} 
The last two examples are generalizations of the Cauchy--Euler equation. 
Differently from the Cauchy--Euler case ($m=0$),  
these generalizations have bounded solutions $\forall {\bf x}_0$,  
provided the conditions discussed in this paper are satisfied. 
 
\vspace{0.5cm}

\begin{center} 
{\bf Acknowledgments} 
\end{center} 
Supported by CNPq/CLAF-150548/2004-4. I appreciate useful discussions  
with Henrique P. Oliveira, Ricardo Mansilla and Eduardo Gallestey.

\end{document}